\documentclass[UTF-8,reqno]{amsart}
\usepackage{enumerate}
\usepackage{amssymb,url,color, booktabs}
\usepackage{mathrsfs}

\numberwithin{equation}{section}

\newcommand{\be}{\begin{eqnarray}}
\newcommand{\ee}{\end{eqnarray}}
\newcommand{\ce}{\begin{eqnarray*}}
\newcommand{\de}{\end{eqnarray*}}
\newtheorem{theorem}{Theorem}[section]
\newtheorem{lemma}[theorem]{Lemma}
\newtheorem{remark}[theorem]{Remark}
\newtheorem{definition}[theorem]{Definition}
\newtheorem{proposition}[theorem]{Proposition}
\newtheorem{Examples}[theorem]{Example}
\newtheorem{corollary}[theorem]{Corollary}

\def\nor{|\mspace{-3mu}|\mspace{-3mu}|}

\def\eps{\varepsilon}

\def\u{\mathbf{u}}

\def\p{\partial}

\def\[{{\Big[}}
\def\]{{\Big]}}
\def\<{{\langle}}
\def\>{{\rangle}}
\def\({{\Big(}}
\def\){{\Big)}}

\def\bx{{\mathbf{x}}}

\def\dif{{\mathord{{\rm d}}}}

\def\no{\nonumber}
\def\={&\!\!=\!\!&}

\def\cF{{\mathcal F}}

\def\cP{{\mathcal P}}

\def\mC{{\mathbb C}}

\def\mE{{\mathbb E}}

\def\mH{{\mathbb H}}
\def\mI{{\mathbb I}}

\def\mL{{\mathbb L}}

\def\mN{{\mathbb N}}

\def\mP{{\mathbb P}}
\def\mQ{{\mathbb Q}}
\def\mR{{\mathbb R}}

\def\mW{{\mathbb W}}

\def\mZ{{\mathbb Z}}

\def\bP{{\mathbf P}}

\def\bE{{\mathbf E}}
\def\1{{\mathbf{1}}}

\def\sF{{\mathscr F}}

\def\sI{{\mathscr I}}

\def\geq{\geqslant}
\def\leq{\leqslant}

\def\div{\mathord{{\rm div}}}

\def\eps{\varepsilon}

\def\u{\mathbf{u}}

\def\p{\partial}

\def\[{{\Big[}}
\def\]{{\Big]}}
\def\<{{\langle}}
\def\>{{\rangle}}
\def\({{\Big(}}
\def\){{\Big)}}

\def\bx{{\mathbf{x}}}

\def\dif{{\mathord{{\rm d}}}}

\def\no{\nonumber}
\def\={&\!\!=\!\!&}
\def\bt{\begin{theorem}}
\def\et{\end{theorem}}
\def\bl{\begin{lemma}}
\def\el{\end{lemma}}
\def\br{\begin{remark}}
\def\er{\end{remark}}
\def\bx{\begin{Examples}}
\def\ex{\end{Examples}}
\def\bd{\begin{definition}}
\def\ed{\end{definition}}
\def\bp{\begin{proposition}}
\def\ep{\end{proposition}}
\def\bc{\begin{corollary}}
\def\ec{\end{corollary}}

\def\wt{\widetilde}

\def\geq{\geqslant}
\def\leq{\leqslant}

\def\div{\mathord{{\rm div}}}

\def\bK{{\mathbf K}}

\def\bP{{\mathbf P}}

\def\<{\langle} \def\>{\rangle}

 \def\beq{\begin{equation}}

\allowdisplaybreaks

\begin{document}

\title{Weak solutions of McKean-Vlasov SDEs with supercritical drifts}
\date{}
\author{Xicheng Zhang}

\address{Xicheng Zhang:
School of Mathematics and Statistics, Wuhan University,
Wuhan, Hubei 430072, P.R.China\\
Email: XichengZhang@gmail.com
 }

\thanks{
This work is supported by NNSFC grant of China (No. 11731009) and the DFG through the CRC 1283 
``Taming uncertainty and profiting from randomness and low regularity in analysis, stochastics and their applications''.
}

\begin{abstract}
Consider the following McKean-Vlasov SDE:
$$
{\mathrm d} X_t=\sqrt{2}\dif  W_t+\int_{\mR^d}K(t,X_t-y)\mu_{X_t}(\dif y){\mathrm d} t,\ \ X_0=x,
$$
where $\mu_{X_t}$ stands for the distribution of $X_t$ and $K(t,x): {\mathbb R}_+\times{\mathbb R}^d\to{\mathbb R}^d$ 
is a time-dependent divergence free vector field. 
Under the assumption $K\in L^q_t(\widetilde L_x^p)$ with $\frac dp+\frac2q<2$, where $\widetilde L^p_x$ stands for the localized $L^p$-space,
we show the existence of weak solutions 
to the above SDE. As an application, we provide a new proof for the existence of weak solutions to 2D-Navier-Stokes equations with 
measure as initial vorticity.


\end{abstract}

\maketitle \rm

\section{Introduction}

Consider the following two dimensional Naver-Stokes equation:
\begin{align}\label{2NS}
\dif\u=\Delta\u+\u\cdot\nabla u+\nabla p,\ \ \div\u=0,
\end{align}
where $\u=(u_1,u_2)$ stands for the velocity field, and $p$ stands for the pressure.
Let $\rho:={\rm curl}\u=\p_1u_2-\p_2 u_1$ be the vorticity of $\u$. It is easy to see that
$$
\p_t\rho=\Delta\rho+\u\cdot\nabla\rho=\Delta\rho+\div(\rho\cdot\u).
$$
Moreover, by the Biot-Savart law we have (cf. \cite{Ma-Be})
$$
\u(t,x)=\int_{\mR^2}K_2(x-y)\rho(t,y)\dif y=:K_2*\rho(t,x),
$$
where
\begin{align}\label{K2}
K_2(x):=\tfrac1{2\pi}\Big(\tfrac{-x_2}{|x|^2},\tfrac{x_1}{|x|^2}\Big).
\end{align}
In other words, $\rho$ solves the following nonlinear integral-differential equation:
\begin{align}\label{2VNS}
\p_t\rho=\Delta\rho+\div(\rho\cdot K_2*\rho).
\end{align}
Notice that the kernel function $K_2$ is of homogeneous of degree $-1$, and
\begin{align}\label{DD1}
\int_{\mR^2}|K_2(x)|^p\dif x=\infty,\ \ p\in[1,\infty].
\end{align}
Suppose that $\rho(0,x)\geq 0$ and $\int_{\mR^2}\rho(0,x)\dif x=1$. 
By the maximum principle and integrating both sides of \eqref{2VNS} with respect to $x$, we obtain
that for any $t>0$,
$$
\rho(t,x)\geq 0,\ \ \int_{\mR^2}\rho(t,x)\dif x=\int_{\mR^2}\rho(0,x)\dif x=1,
$$
which means that $\rho(t,\cdot)_{t\geq 0}$ is a family of probability measures.
By the superposition principle \cite{Tr}, there would be a weak solution to the following McKean-Vlasov SDEs:
\begin{align}\label{SDE0}
\dif X_t=\left[\int_{\mR^2}K_2(X_t-y)\rho(t,y)\dif y\right]\dif t+\sqrt{2}\dif W_t,\ X_0=x,
\end{align}
where $W$ is a two-dimensional standard Brownian motion and $\rho(t,\cdot)$ is the distributional density of $X_t$. 
On the other hand, if $X_t$ solves the nonlinear SDE \eqref{SDE0}, by It\^o's formula, 
the law of $X_t$ will solve the nonlinear Fokker-Planck equation \eqref{2VNS} in the distributional sense. It should be noticed that
when the initial vorticity is a finite Radon measure, the existence of solutions to PDE \eqref{2VNS} was obtained by 
Giga, Miyakawa and Osada in \cite{Gi-Mi-Os} (see also Cottet's work \cite{Co}), and the uniqueness was proven 
by Gallagher and Gallay in \cite{Ga-Ga}.
 However, due to the non-integrability of $K_2$ (see \eqref{DD1}), 
it does not immedately imply the existence of weak solutions to distributional dependent SDE \eqref{SDE0} 
by superposition principle because the following condition in \cite{Tr} is not known to hold for the solution obtained in \cite{Gi-Mi-Os},
$$
\int^T_0\left|\int_{\mR^2}\left[\int_{\mR^2}K_2(x-y)\rho(t,y)\dif y\right]\rho(t,x)\dif x\right|\dif t<\infty.
$$

In this paper we are concerning with the following McKean-Vlasov SDE in $\mR^d$:
\begin{align}\label{SDE}
\dif X_t=\left[\int_{\mR^d}K(t, X_t,y)\mu_{X_t}(\dif y)\right]\dif t+\sqrt{2}\dif W_t,
\end{align}
where $K:\mR_+\times\mR^d\times\mR^d\to\mR^d$ is a measurable vector-valued function and $\mu_{X_t}$ is the law of $X_t$. 
For any $\alpha\in[0,2)$, we introduce the following index set:
$$
\sI_\alpha:=\Big\{(p,q)\in(1,\infty)^2,\ \tfrac{d}{p}+\tfrac{2}q<2-\alpha\Big\}.
$$
Suppose that for some $(p,q)\in\sI_1$,
\begin{align}\label{Sub}
|K(t,x,y)|\leq h(t,x-y),\ \ h\in  L^q_t(\widetilde L^p_x):=\cap_{T>0}L^q([0,T];\widetilde L^p),
\end{align}
where $\widetilde L^p$ is the localized $L^p$-space in $\mR^d$ (see \eqref{Ck} below).
Under \eqref{Sub}, R\"ockner and the present author \cite{Ro-Zh} showed the strong well-posedness to the above SDE.
The integrability condition \eqref{Sub} for $\frac dp+\frac2q<1$ is usually called subcritical case in the literature; while $\frac dp+\frac2q=1$ 
and $\frac dp+\frac2q>1$
correspond to the critical and supercritical cases, respectively.
Notice that the kernel function $K_2$ given in \eqref{K2} belongs to the supercritical regime since
$$
\int_{|x|<1}|K_2(x)|^p\dif x<\infty,\ \ p\in[1,2),\ \ \int_{|x|<1}|K_2(x)|^2\dif x=\infty.
$$

For $\beta\geq 0$, let $\cP_\beta(\mR^d)$ be the set of all probability measures on $\mR^d$ with finite $\beta$-order moment.
The aim of this paper is to show the following existence result.
\bt\label{TT1}
We suppose that  in the distributional sense,
\begin{align}\label{AA1}
\div K(t,\cdot,y)\leq 0,
\end{align}
and for some $(p,q)\in\sI_0$,
\begin{align}\label{AA2}
|K(t,x,y)|\leq h(t,x-y),\ \ h\in  L^q_t(\widetilde L^p_x).
\end{align}
Let $\beta\in[0,2/(\frac dp+\frac 2q))$. For any $\nu_0\in\cP_\beta(\mR^d)$, there exists at least one weak solution to SDE \eqref{SDE} with initial distribution $\nu_0$.
More precisely, there are stochastic basis $(\Omega,\sF,(\sF_t)_{t\geq 0}, \bP)$ and two $\sF_t$-adapted processes $(X,W)$ defined on it such that
\begin{enumerate}[(i)]
\item $\bP\circ X^{-1}_0=\nu_0$ and $W$ is a $d$-dimensional standard $\sF_t$-Brownian motion.
\item It holds that for all $t\geq 0$,
$$
X_t=X_0+\int^t_0\!\!\!\int_{\mR^d}K(s, X_s,y)\mu_{X_s}(\dif y)\dif s+\sqrt{2}W_t,\ \bP-a.s.,
$$
where $\mu_{X_s}$ is the law of $X_s$.
\end{enumerate}
Moreover, we have the following conclusions:
\begin{enumerate}[(i)]
\item For any $T>0$, there is a constant $C>0$ such that
\begin{align}\label{Mo9}
\bE\left(\sup_{t\in[0,T]}|X_t|^\beta\right)\leq C(\bE|X_0|^\beta+1).
\end{align}

\item For Lebesgue almost all $t>0$, $X_t$ admits a density $\rho(t,\cdot)$ with the regularity
\begin{align}\label{Reg}
\rho\in \cap_{T>0}\mH^{\alpha,p}_q(T),\ \alpha\in[0,1],\ p,q\in(1,\infty),\ \tfrac dp+\tfrac 2q>d+\alpha,
\end{align}
where $\mH^{\alpha,p}_q(T):=L^q([0,T]; H^{\alpha,p})$ and $H^{\alpha,p}$ is the usual Bessel potential space.
\end{enumerate}
\et

Recently, there are great interests to study the McKean-Vlasov or distributional dependent SDEs since it appears in the studies of 
propagation of chaos \cite{Sz}, mean-field games (cf. \cite{Ca-De1}) and nonlinear integral-partial differential equations \cite{Mc, Fu}.
When $K$ is bounded measurable, the existence and uniqueness of weak solutions to SDE \eqref{SDE} 
was proved by Li and Min \cite{Li-Mi} (see also \cite{Mi-Ve} for the strong well-posedness of SDE \eqref{SDE}). 
As mentioned above, when $K$ is singular and belongs to the subcritical regime, the strong existence and uniqueness
was shown in \cite{Ro-Zh} recently (see also \cite{Hu-Wa}).
We also mention that Jabin and Wang  \cite{Ja-Wa} showed the propagation of chaos for singular 
kernel $K_2$ above by purely analytic method.
While, the existence of particle trajectories is not provided therein.
Here, an open question is the uniqueness of weak solutions in the supercritical case. 
This is even not known for linear SDEs with supercritical drifts (cf. \cite{Zh-Zh2}).

\medskip

As a simple application of Theorem \ref{TT1}, we have the following corollary.
\bc
Consider the vorticity form \eqref{2VNS} of $2D$-Navier-Stokes equations. Let $\beta\in[0,2)$.
For any $\mu(0)\in\cP_\beta(\mR^2)$, there exists a continuous curve $t\mapsto\mu(t)\in\cP_\beta(\mR^2)$
such that for all $t\geq 0$ and $f\in C^\infty_b(\mR^d)$,
$$
\mu_t(f)=\mu_0(f)+\int^t_0\mu_s(\Delta f)\dif s+\int^t_0\left[\int_{\mR^2}\!\!\int_{\mR^2}K_2(x-y)\cdot\nabla f(x)\mu_s(\dif y)\mu_s(\dif x)\right]\dif s.
$$
Moreover, for Lebesgue-almost all $t>0$, $\mu(t,\dif x)=\rho(t,x)\dif x$, where $\rho$ satisfies \eqref{Reg}, 
and for any $T>0$ and some $C>0$,
\begin{align}\label{Mom}
{\rm ess.}\sup_{t\in[0,T]}\int_{\mR^2}|x|^\beta\rho(t,x)\dif x\leq C\left(\int_{\mR^2}|x|^\beta\mu_0(\dif x)+1\right).
\end{align}
\ec

\br\rm
Compared with \cite{Gi-Mi-Os}, the new point here is that the moment estimate \eqref{Mom} is obtained, 
which provides the decay estimate of the vorticity as $|x|\to\infty$. Note that the uniqueness is proven in \cite{Ga-Ga}, 
which strongly depends on the structure of $K_2$.
\er

This paper is organized as following: in Section 2, we prepare necessary spaces and 
some well-known results about the maximum principle for the associated PDE.
In Section 3, through mollifying the kernel function $K$, we show our main result by  weak convergence method.

\section{Preliminaries}

We first introduce the following spaces and notations for later use.
For $(\alpha,p)\in\mR_+\times[1,\infty]$,  the  Bessel potential space $H^{\alpha,p}$ is defined by
$$
H^{\alpha,p}:=\big\{f\in L^1_{loc}(\mR^d): \|f\|_{\alpha,p}:=\|(\mI-\Delta)^{\alpha/2}f\|_p<\infty\big\},
$$
where $\|\cdot\|_p$ is the usual $L^p$-norm, and $(\mI-\Delta)^{\alpha/2}f$ is defined by Fourier's  transform
$$
(\mI-\Delta)^{\alpha/2}f:=\cF^{-1}\big((1+|\cdot|^2)^{\alpha/2}\cF f\big).
$$
For $T>0$, $p,q\in[1,\infty]$ and $\alpha\in\mR_+$, we introduce the following spaces of space-time functions,
$$
\mL^p_q(T):=L^q\big([0,T];L^p\big),\ \  \mH^{\alpha,p}_q(T):=L^q\big([0,T];H^{\alpha,p}\big).
$$
Let $\chi\in C^\infty_c(\mR^d)$ be a smooth function with $\chi(x)=1$ for $|x|\leq 1$ and $\chi(x)=0$ for $|x|>2$.
For $r>0$ and $z\in\mR^d$, define
\begin{align}\label{CHI}
\chi^z_r(x):=\chi((x-z)/r).
\end{align}
Fix $r>0$. We introduce the following localized $H^{\alpha,p}$-space:
\begin{align}\label{Ck}
\widetilde H^{\alpha,p}:=\Big\{f: \nor f\nor_{\alpha,p}:=\sup_z\|f\chi^z_r\|_{\alpha,p}<\infty\Big\},
\end{align}
and  the localized space-time function space 
\begin{align}\label{GG1}
\widetilde\mH^{\alpha,p}_q(T):=\Big\{f: \nor f\nor_{\widetilde\mH^{\alpha,p}_q(T)}:=\sup_{z\in\mR^d}\|\chi^z_r f\|_{\mH^{\alpha,p}_q(T)}<\infty.\Big\}
\end{align}
For simplicity we shall write
$$
\widetilde\mH^{\alpha,p}_q:=\cap_{T>0}\widetilde\mH^{\alpha,p}_q(T),\ \ \widetilde\mL^{p}_q:=\widetilde\mH^{0,p}_q,\ \ \wt\mL^p(T):=\wt\mL^p_p(T).
$$
It should be noticed that
$$
L^q([0,T];\widetilde L^p)\subset\widetilde\mL^p_q(T).
$$

The following lemma lists some easy properties of $\widetilde\mH^{\alpha,p}_q$ (see \cite{Zh-Zh2}).
\bp\label{Pr41}
Let $p,q\in(1,\infty)$, $\alpha\in\mR_+$ and $T>0$.
\begin{enumerate}[{\rm(i)}]
\item For $r\not=r'>0$, there is a $C=C(d,\alpha,r,r',p,q)\geq 1$ such that% for all $f\in\widetilde \mH^{\alpha,p}_q$,
\begin{align}\label{GT1}
C^{-1}\sup_{z}\|f\chi^{z}_{r'}\|_{\mH^{\alpha,p}_q(T)}\leq \sup_{z}\|f\chi^{z}_r\|_{\mH^{\alpha,p}_q(T)}\leq C \sup_{z}\|f\chi^{z}_{r'}\|_{\mH^{\alpha,p}_q(T)}.
\end{align}
In other words, the definition of $\widetilde \mH^{\alpha,p}_q$ does not depend on the choice of $r$.
\item Let $(\rho_\eps)_{\eps\in(0,1)}$ be a family of mollifiers in $\mR^d$. 
For any $f\in\widetilde \mH^{\alpha,p}_q$ and $T>0$, it holds that 
$$
f_\eps(t,x):=f(t,\cdot)*\rho_\eps(x)\in L^q_{loc}(\mR; C^\infty_b(\mR^d)),
$$
 and for some $C=C(d,\alpha,p,q)>0$,
\begin{align}\label{GT2}
\nor f_\eps\nor_{\widetilde\mH^{\alpha,p}_q(T)}\leq C\nor f\nor_{\widetilde\mH^{\alpha,p}_q(T)},\ \forall \eps\in(0,1),
\end{align}
and for any $\varphi\in C^\infty_c(\mR^{d})$,
\begin{align}\label{GT3}
\lim_{\eps\to0}\|(f_\eps-f)\varphi\|_{\mH^{\alpha,p}_q(T)}=0.
\end{align}
\end{enumerate}
\ep
We introduce the following notion about Krylov's estimate.
\bd
Let $p,q\in(1,\infty)$ and $T,\kappa>0$. We say a stochastic process $X$ 
satisfies Krylov's estimate with index $p,q$ and constant $\kappa$ if for any 
$f\in \widetilde\mL^p_q(T)$,
\begin{align}\label{Kry}
\bE\left(\int^T_0f(t,X_t)\dif t\right)\leq \kappa\nor f\nor_{\widetilde\mL^p_q(T)}.
\end{align}
The set of all such $X$ will be denoted by ${\bf K}^{p,q}_{T,\kappa}$.
\ed
\br\label{Re9}\rm
By Krylov's estimate \eqref{Kry}, there is a density function $\rho^X\in\mL^{r}_{s}(T)$ with $r=\frac{p}{p-1}$ and $s=\frac{q}{q-1}$ so that
$$
\int^T_0\!\!\!\int_{\mR^d}f(t,x)\rho^X_t(x)\dif x\dif t=\bE\left(\int^T_0f(t,X_t)\dif t\right)
\leq \kappa\nor f\nor_{\widetilde\mL^{p}_{q}(T)}\leq \kappa\|f\|_{\mL^{p}_{q}(T)}.
$$
\iffalse
By Proposition \ref{Pr41} (v), we further have
$$
\nor\rho^X\nor^*_{\widetilde\mL^{r}_{s}(T)}:=\sum_{z\in\mZ^d}\|\1_{Q_z}\rho^X\|_{\mL^{r}_{s}(T)}\leq \kappa,
$$
where $Q_z$ is defined by \eqref{QQZ}.
\fi
\er
\iffalse%%%%%%%%%%%%%%%%%%%%%%%%%%%%%%%
\br
Let $X\in{\bf K}^{p,q}_{T,\kappa}$. By the interpolation theorem, for any $p'\in[p,\infty)$ and $q'\in(1,q]$ with
$\tfrac{p'}{q'}-p'=\tfrac{p}{q}-p$, we have $X\in{\bf K}^{p',q'}_{T,\kappa'}$ for some $\kappa'>0$.
In fact, note that by the interpolation theorem (see \cite[Theorem 5.1.2]{Be-Lo}),
$$
(\widetilde\mL^\infty_1(T), \widetilde\mL^p_q(T))_{[\theta]}=\widetilde\mL^{p'}_{q'}(T),
$$
where $\theta\in(0,1)$, $\frac{1}{q'}=1-\theta+\frac{\theta}{q}$ and $\frac{1}{p'}=\frac{\theta}{p}$, $(\cdot,\cdot)_{[\theta]}$ stands for the complex interpolation.
On the other hand, since obviously,
$$
\bE\left(\int^T_0f_t(X_t)\dif t\right)\leq\nor f\nor_{\widetilde\mL^\infty_1(T)}.
$$
%$$\nor f\nor_{\widetilde\mL^{p'}_{q'}(T)}\leq \nor f\nor^{p/p'}_{\widetilde\mL^{p}_{q}(T)}\nor f\nor^{1-p/p'}_{\widetilde\mL^{\infty}_{1}(T)}.$$
By the interpolation theorem, we get $X\in{\bf K}^{p',q'}_{T,\kappa'}$.
\er
\fi %%%%%%%%%%%%%%%%%%%%
For a space-time function $f(t,x,y):\mR_+\times\mR^d\times\mR^d\to\mR$ and $p_1,p_2,q_0\in[1,\infty]$, we also introduce the norm
\begin{align}\label{De9}
\nor f\nor_{\widetilde\mL^{p_1,p_2}_{q_0}(T)}:=\sup_{z,z'\in\mR^d}\left(\int^T_0\left(\int_{\mR^d}\1_{B^{z'}_1(y)}\|\1_{B^z_1}f(t,\cdot,y)\|^{p_2}_{p_1}\dif y\right)^{\frac{q_0}{p_2}}\right)^{\frac{1}{q_0}},
\end{align}
where for $z\in\mR^d$ and $r>0$,
$$
B^z_r:=\{x\in\mR^d: |x-z|<r\},\ \ B_r:=B^0_r.
$$
The following lemma will be used to take the limits in the proof of the existence of weak solutions (see \cite[Lemma 2.6]{Ro-Zh}). 
\bl\label{Le27}
Let $p_1,p_2, q_0,q_1,q_2\in(1,\infty)$ with $\frac{1}{q_1}+\frac{1}{q_2}=1+\frac{1}{q_0}$ and $T,\kappa_1,\kappa_2>0$. 
Let $X\in\bK^{p_1,q_1}_{T,\kappa_1}$ and $Y\in\bK^{p_2,q_2}_{T,\kappa_2}$ be two independent processes.
Then for any $f(t,x,y)\in\widetilde\mL^{p_1,p_2}_{q_0}(T)$,
\begin{align}\label{GQ2}
\bE\left(\int^T_0f(t,X_t,Y_t)\dif t\right)&\leq \kappa_1\kappa_2\nor f\nor_{\widetilde\mL^{p_1,p_2}_{q_0}(T)}.
\end{align}
\el

Consider the following backward PDE:
\begin{align}\label{PDE0}
\p_t u+\Delta u+b\cdot\nabla u=f,\ \ u(T)\equiv0.
\end{align}
The following maximum principle was proven in \cite[Theorem 2.2]{Zh-Zh2}.
\bt\label{pde}
Let $T>0$. Suppose that $b\in C^\infty_b([0,T]\times\mR^d))$ satisfies $\div b\leq 0$ and for some $(p,q)\in\sI_0$ and $\kappa>0$,
$$
\nor b\nor_{\widetilde\mL^p_q(T)}\leq \kappa. 
$$
Let $\alpha\in[0,1]$ and $f\in C^\infty_0(\mR^{d+1})$. 
For any $(\bar p,\bar q)\in\sI_\alpha$, there is a constant $C>0$ only depending on $T,d,p,q,\alpha,\bar p,\bar q,\kappa$ 
such that for any smooth solution $u$ of PDE \eqref{PDE0},
\begin{align}\label{Max}
\|u\|_{L^\infty([0,T]\times\mR^d)}\leq C\nor f\nor_{\widetilde\mH^{-\alpha,\bar p}_{\bar q}(T)}.
\end{align}
\et

\section{Proof of Theorem \ref{TT1}}

Suppose that $K(t,x,y)$ satisfies \eqref{AA1} and \eqref{AA2}.
Let $(\varrho^{\rm d}_n)_{n\in\mN}$ be a family of mollifiers in $\mR^d$ with compact supports in unit ball $B_1$.
Define
\begin{align}\label{D8}
K_n(t,x,y):=\int^\infty_0\!\!\!\int_{\mR^{2d}}K(t',x',y')\varrho^{\rm 1}_n(t-t')\varrho^{\rm d}_n(x-x')\varrho^{\rm d}_n(y-y')\dif t'\dif x'\dif y',
\end{align}
and for a probability measure $\mu$,
$$
b_n(t,x,\mu):=\int_{\mR^d}K_n(t,x,y)\mu(\dif y).
$$
By (ii) of Proposition \ref{Pr41} and \eqref{AA2}, one sees that for each $T>0$ and $j=0,1,\cdots,$
\begin{align}\label{KK0}
\kappa^j_n:=\|\nabla^j_x K_n\|_{\mL^\infty(T)}+\|\nabla^j_y K_n\|_{\mL^\infty(T)}<\infty
\end{align}
and
$$
\div b_n(t,\cdot,\mu)\leq 0.
$$
Hence, for any $T>0$ and some $C=C(\kappa^1_n)>0$,
$$
|b_n(t,x,\mu)-b_n(t,\bar x,\mu)|\leq C |x-\bar x|,
$$
and for any two random variables $X,Y$,
$$
|b_n(t,x,\mu_X)-b_n(t,x,\mu_Y)|\leq C \mE|X-Y|.
$$
Thus, there is a unique strong solution $X^n_t$ to the following McKean-Vlasov SDE:
\begin{align}\label{SDE1}
\dif X^n_t=b_n(t, X^n_t,\mu_{X^n_t})\dif t+\sqrt{2}\dif W_t,\ X^n_0\stackrel{(d)}{=}\nu_0,
\end{align}
where $W$ is a $d$-dimensional Brownian motion on some probability space $(\Omega,\sF,\mP)$.

We first show the following key Krylov estimate.
\bl
Let $\alpha\in[0,1]$ and $(\bar p,\bar q)\in\sI_\alpha$. For any $T>0$, there is a constant $C>0$ such that for each $n\in\mN$,
\begin{align}\label{Kry0}
\mE\left(\int^T_0f(t,X^n_t)\dif t\right)\leq C\nor f\nor_{\widetilde\mH^{-\alpha,\bar p}_{\bar q}(T)}.
\end{align}
In particular, $X^n\in\bK^{\bar p,\bar q}_{T,\kappa}$ for any $(\bar p,\bar q)\in\sI_0$.
\el
\begin{proof}
Without loss of generality, we assume that $f$ is smooth.
By the properties of convolution and $|K(t,x,y)|\leq h(t,x-y)$, one sees that
$$
\nor b_n(t,\cdot,\mu)\nor_p\leq \int_{\mR^d}\nor K_n(t,\cdot,y)\nor_p \mu(\dif y)\leq \int^\infty_0\nor h(t',\cdot)\nor_p\varrho^{\rm 1}_n(t-t')\dif t'.
$$
Hence,
\begin{align}\label{AA4}
\int^T_0\nor b_n(t,\cdot,\mu_{X^n_t})\nor_p^q\dif t\leq \int^T_0\nor h(t,\cdot)\nor_p^q\dif t.
\end{align}
Now, consider the following backward PDE:
$$
\p u_n+\Delta u_n+b_n(t,\cdot,\mu_{X^n_t})\cdot\nabla u_n=f,\ \ u_n(T)=0.
$$
Since by \eqref{KK0}, for any $j\in\mN$, 
$$
\sup_{t\in[0,T]}\|\nabla^j b_n(t,\cdot,\mu_{X^n_t})\|_\infty<\infty,
$$
the above PDE admits a unique smooth solution $u_n$ with the regularities
$$
\p_t u_n, \nabla^2 u_n\in L^\infty([0,T]\times\mR^d).
$$
Moreover, by Theorem \ref{pde}, there is a constant $C>0$ such that for all $n\in\mN$,
$$
\|u_n\|_{L^\infty([0,T]\times\mR^d)}\leq C\nor f\nor_{\widetilde\mH^{-\alpha,\bar p}_{\bar q}(T)}.
$$
Now by It\^o's formula, we have
$$
\mE u_n(T, X^n_T)=\mE u_n(0, X^n_0)+\mE\int^T_0f(t, X^n_t)\dif t,
$$
which implies that
$$
\mE\int^T_0 f(s, X^n_s)\dif s\leq \|u_n(0, \cdot)\|_\infty\leq C\nor f\nor_{\widetilde\mH^{-\alpha,\bar p}_{\bar q}(T)}.
$$
The proof is complete.
\end{proof}

Let $\mC$ be the space of all continuous functions from $\mR_+$ to $\mR^d$, which is endowed with the locally uniform convergence topology so that $\mC$ is a Polish space. We also use the following convention below: 
The letter $C$ with or without subscripts will denote a constant whose value may change in different places.

\bl\label{Le32}
Let $\mP_n$ be the law of $X^n_\cdot$ in $\mC$. Then $(\mP_n)_{n\in\mN}$ is tight. Moreover, for any $\beta\in[0,2/(\frac dp+\frac 2q))$ and $T>0$, 
there is a constant $C>0$ such that for all $n\in\mN$,
\begin{align}\label{Mo}
\mE\left(\sup_{t\in[0,T]}|X^n_t|^\beta\right)\leq C(\mE|X^n_0|^\beta+1).
\end{align}
\el
\begin{proof}
Let $T>0$ and $\tau\leq T$ be any bounded stopping time. 
For any $\delta>0$ and $\gamma\in(1,2/(\frac dp+\frac 2q))$, by H\"older's inequality and \eqref{Kry0} with $\alpha=0$, we have
\begin{align}
\mE\left(\int^{\tau+\delta}_\tau |b_n(t, X^n_t,\mu_{X^n_t})|\dif t\right)^\gamma
&\leq \delta^{\gamma-1}\mE\left(\int^{\tau+\delta}_\tau |b_n(t, X^n_t,\mu_{X^n_t})|^\gamma\dif t\right)\no\\
&\leq \delta^{\gamma-1}\mE\left(\int^{T+\delta}_0 |b_n(t, X^n_t,\mu_{X^n_t})|^\gamma\dif t\right)\no\\
&\leq C\delta^{\gamma-1}\left[\int^{T+\delta}_0 \nor |b_n(t,\cdot,\mu_{X^n_t})|^\gamma\nor_{p/\gamma}^{q/\gamma}\dif t\right]^{\frac\gamma{q}}\no\\
&= C\delta^{\gamma-1}\left[\int^{T+\delta}_0 \nor b_n(t,\cdot,\mu_{X^n_t})\nor_{p}^{q}\dif t\right]^{\frac\gamma{q}}\no\\
&\stackrel{\eqref{AA4}}{\leq} C\delta^{\gamma-1}\left[\int^{T+\delta}_0 \nor h(t,\cdot)\nor_{p}^{q}\dif t\right]^{\frac\gamma{q}},\label{Mo1}
\end{align}
where the constant $C$ does not depend on $n$. Moreover, by Burkholder's inequality, it is easy to see that
$$
\mE\left(\sup_{s\in[0,\delta]}|W_{\tau+s}-W_\tau|^\gamma\right)\leq C\delta^{\gamma/2}.
$$
Hence, for any $\theta\in(0,1)$, by \cite[Lemma 2.7]{Zh-Zh1} we have
$$
\sup_n\mE\left(\sup_{s,t\in[0,T],s\not=t}|X^n_t-X^n_s|^{\theta\gamma}\right)\leq C\delta^{\theta(\gamma-1)}.
$$
The tightness of $(\mP_n)_{n\in\mN}$ now follows by \cite[Theorem 1.3.2]{St-Va}. 
Finally, the moment estimate \eqref{Mo} follows by \eqref{SDE1} and \eqref{Mo1}.
\end{proof}

Now we can give
\begin{proof}[Proof of Theorem \ref{TT1}]
Let $\mP_n$ be the law of $X^n_\cdot$ in $\mC$ and $\mW$ the law of Brownian motion in $\mC$. Consider the product probability measure $\mQ_n:=\mP_n\times\mP_n\times\mW$.
By Lemma \ref{Le32}, one sees that $(\mQ_n)_{n\in\mN}$ is tight in $\mC\times\mC\times\mC$.
Let $\mQ$ be any accumulation point. 
Without loss of generality, we assume that $\mQ_n$ weakly converges to some probability measure $\mQ$.
By Skorokhod's representation theorem,  there are a probability space $(\tilde\Omega,\tilde\sF,\tilde\bP)$ and 
random variables $(\tilde X^n , \tilde Y^n, \tilde W^n)$ and $(\tilde X,\tilde Y,\tilde W)$ defined on it such that
\begin{align}\label{FD4}
(\tilde X^n, \tilde Y^n, \tilde W^n)\to (\tilde X,\tilde Y, \tilde W),\ \ \tilde\bP-a.s.
\end{align}
and
\begin{align}\label{FD5}
\tilde\bP\circ(\tilde X^n,\tilde Y^n,  \tilde W^n)^{-1}=\mQ_n,\quad
\tilde\bP\circ(\tilde X, \tilde Y, \tilde W)^{-1}=\mQ.
\end{align}
Define $\tilde\sF^n _t:=\sigma(\tilde W^n _s, \tilde X^n_s;s\leq t)$. 
We note that
\begin{align*}
&\mP(W _t-W _s\in\cdot |\sF _s)=\mP(W _t-W _s\in\cdot)\\
&\Rightarrow
\tilde\bP(\tilde W^n _t-\tilde W^n _s\in\cdot |\tilde \sF^n _s)=\tilde\bP(\tilde W^n _t-\tilde W^n _s\in\cdot).
\end{align*}
In other words, $\tilde W^n $ is an $\tilde\sF_t^n $-Brownian motion. Thus, by \eqref{SDE1} and \eqref{FD5} we have
\begin{align}\label{DD3}
\tilde X^n _t=\tilde X^n_0
+\int^t_0b_n(s,\tilde X^n_s,\mu_{\tilde X^n_s})\dif s+\sqrt{2}\tilde W^n _t.
\end{align}
To show the existence of a solution, the key point is to show that
\begin{align}\label{DD4}
\int^t_0b_n(s,\tilde X^n_s,\mu_{\tilde X^n_s})\dif s\to \int^t_0b(s,\tilde X_s,\mu_{\tilde X_s})\dif s\mbox{ in probability as $n\to\infty$},
\end{align}
where $b(s,x,\mu):=\int_{\mR^d}K(s,x,y)\mu(\dif y)$.
After showing this limit, we can take limits for both sides of \eqref{DD3} to obtain the existence of a solution, i.e.
$$
\tilde X _t=\tilde X_0
+\int^t_0b(s,\tilde X_s,\mu_{\tilde X_s})\dif s+\sqrt{2}\tilde W _t.
$$
Since $\tilde X^n$ and $\tilde Y^n$ are independent by \eqref{FD5}, to prove \eqref{DD4}, it suffices to show that 
$$
\int^t_0K_n(s,\tilde X^n_s,\tilde Y^n_s)\dif s\to \int^t_0K(s,\tilde X_s,\tilde Y_s)\dif s\mbox{ in probability as $n\to\infty$}.
$$
The above limit will be a consequence of the following two limits: for each $m\in\mN$,
\begin{align}\label{Lim1}
\lim_{n\to\infty}\int^t_0|K_m(s, \tilde X^n_s, \tilde Y^n_s)-K_m(s,\tilde X_s,\tilde Y_s)|\dif s=0, \  \bP-a.s.
\end{align}
and
\begin{align}\label{Lim2}
\lim_{m\to\infty}\sup_n\bE\int^t_0|K_m(s, \tilde X^n_s, \tilde Y^n_s)-K(s, \tilde X^n_s, \tilde Y^n_s)|\dif s=0.
\end{align}
Below we drop the tilde for simplicity. For fixed $m\in\mN$,
since $K_m$ is bounded and $(x,y)\mapsto K_m(s,x,y)$ is continuous, it follows by the dominated convergence theorem and \eqref{FD4} 
that the limit \eqref{Lim1} holds. For limit \eqref{Lim2}, we write
$$
\bE\int^t_0|K_m(s, X^n_s, Y^n_s)-K(s, X^n_s, Y^n_s)|\dif s=I^{(1)}_{n,m}(R)+I^{(2)}_{n,m}(R),
$$
where
\begin{align*}
I^{(1)}_{n,m}(R)&:=\bE\int^t_0\1_{\{|X^n_s|\leq R\}\cap\{\ |Y^n_s|\leq R\}}|K_m(s, X^n_s, Y^n_s)-K(s, X^n_s, Y^n_s)|\dif s,\\
I^{(2)}_{n,m}(R)&:=\bE\int^t_0\1_{\{|X^n_s|>R\}\cup\{|Y^n_s|>R\}}|K_m(s, X^n_s, Y^n_s)-K(s, X^n_s, Y^n_s)|\dif s.
\end{align*}
For $I^{(1)}_{n,m}(R)$, since $(p,q)\in\sI_0$, one can choose $\gamma>1$ such that
$$
\tfrac{d}{p}+\tfrac{2\gamma}{q}<2.
$$
Thus, by Lemma \ref{Le27} with $p_1=p$, $q_1=\frac q\gamma$, $p_2>\frac{qd}{2(\gamma-1)}$, $q_2=\frac{q}{q+1-\gamma}$ and $q_0=q$,
\begin{align}\label{G7}
I^{(1)}_{n,m}(R)\leq C\nor\1_{B_R\times B_R}(K_m-K)\nor_{\widetilde\mL^{p,p_2}_q(t)},
\end{align}
where $C$ is independent of $n, m$. Recalling the definition \eqref{De9}, we further have
$$
I^{(1)}_{n,m}(R)\leq C\left(\int^t_0\left(\int_{B_{R+1}}\|\1_{B_{R+1}}|K_m(s,\cdot,y)-K(s,\cdot,y)\|^{p_2}_p\dif y\right)^{\frac{q}{p_2}}\dif s\right)^{1/q}.
$$
Since $|K(s,x,y)|\leq h(s,x-y)$, by \eqref{D8} we have
\begin{align}\label{A7}
&\int^t_0\left(\sup_m\int_{B_{R+1}}\|\1_{B_{R+1}}|K_m(s,\cdot,y)\|^{p_2}_p\dif y\right)^{\frac{q}{p_2}}\dif s\no\\
&\qquad\leq\int^t_0\left(\int_{B_{R+2}}\|\1_{B_{R+2}}|h(s,\cdot-y)\|^{p_2}_p\dif y\right)^{\frac{q}{p_2}}\dif s\no\\
&\qquad\leq C_R\int^t_0\left(\int_{B_{2(R+2)}}|h(s,x)|^p\dif x\right)^{\frac qp}\dif s<\infty.
\end{align}
Hence, by the dominated convergence theorem, for each $R>0$,
\begin{align}\label{G9}
\lim_{m\to\infty}\sup_nI^{(1)}_{n,m}(R)=0.
\end{align}
For $I^{(2)}_{n,m}(R)$, letting $\alpha\in(1, 2/(\frac{d}{p}+\frac{2}{q}))$, by H\"older's inequality and Lemma \ref{Le27} again,
\begin{align}
I^{(2)}_{n,m}(R)&\leq \int^t_0(\bE|K_m(s, X^n_s, Y^n_s)-K(s, X^n_s, Y^n_s)|^\alpha)^{\frac1\alpha}\no\\
&\qquad\times \bP(\{|X^n_s|>R\}\cup\{|Y^n_s|>R\})^{1-\frac1\alpha}\dif s\no\\
&\leq \left(\int^t_0\bE|K_m(s, X^n_s, Y^n_s)-K(s, X^n_s, Y^n_s)|^\alpha\dif s\right)^{\frac1\alpha}\no\\
&\quad\times \sup_{s\in[0,t]}\bP(\{|X^n_s|>R\}\cup\{|Y^n_s|>R\})^{1-\frac1\alpha}\no\\
&\leq C \left(\nor K_m\nor_{\widetilde\mL^{p,p_2}_q(t)}+\nor K\nor_{\widetilde\mL^{p,p_2}_q(t)}\right)
\sup_{s\in[0,t]}(2\bP\{|X^n_s|>R\})^{1-\frac1\alpha},\label{A8}
\end{align}
where $C$ is independent of $n,m$ and $R$, and $p_2$ is chosen being large enough as in \eqref{G7}.
As in \eqref{A7}, we have
\begin{align}
\nor K_m\nor^q_{\widetilde\mL^{p,p_2}_q(T)}
&= \sup_{z,z'\in\mR^d}\int^T_0\left(\int_{\mR^d}\1_{B^{z'}_1}(y)\|\1_{B^z_1}K_m(t,\cdot,y)\|^{p_2}_{p}\dif y\right)^{\frac{q}{p_2}}\dif t\no\\
&\leq  \sup_{z,z'\in\mR^d}\int^T_0\left(\int_{\mR^d}\1_{B^{z'}_2}(y)\|\1_{B^z_2}h(t,\cdot-y)\|^{p_2}_{p}\dif y\right)^{\frac{q}{p_2}}\dif t\no\\
&\leq C\sup_{z,z'\in\mR^d}\int^T_0\sup_{|y-z'|\leq 2}\|\1_{B^z_2}h(t,\cdot-y)\|_{p}^q\dif t\no\\
&\leq C \int^T_0\nor h(t)\nor_{p}^q\dif t<\infty.\label{A5}
\end{align}
Moreover, by \eqref{DD3} and Chebyschev's inequality, we have
\begin{align}
\sup_{s\in[0,t]}\bP\{|X^n_s|>R\}&\leq\bP\{|X^n_0|>\tfrac R3\}+\bP\left\{\sup_{s\in[0,t]}\sqrt{2}|W_s|>\tfrac R3\right\}\no\\
&\quad+\bP\left\{\int^t_0|b_n(s, X^n_s,\mu_{X^n_s})|\dif s>\tfrac R3\right\}\no\\
&\leq\nu\Big\{|x|>\tfrac R3\Big\}+\frac{C}{R}+\frac{3}R\bE\int^t_0|b_n(s, X^n_s,\mu_{X^n_s})|\dif s\no\\
&\leq\nu\Big\{|x|>\tfrac R3\Big\}+\frac{C}{R},\label{A4}
\end{align}
where the constant $C$ is independent of $R$ and $n$, and the last step is due to \eqref{Mo1}.
Combining \eqref{A8}, \eqref{A5} and \eqref{A4}, we obtain
$$
\lim_{R\to\infty}\sup_{n,m}I^{(2)}_{n,m}(R)=0,
$$
which together with \eqref{G9} yields \eqref{Lim2}.
Moreover, the estimate \eqref{Mo9} follows by \eqref{Mo}, and the regularity estimate \eqref{Reg} follows by \eqref{Kry0} and Remark \ref{Re9}.
The proof is thus complete.
\end{proof}


\begin{thebibliography}{999}

\bibitem{Ca-De1}Carmona R. and Delarue F.: 
{\it Probabilistic theory of mean field games with applications. II. Mean field games with common noise and master 
equations.} Probability Theory and Stochastic Modeling, 84. Springer, 2018.

\bibitem{Co} Cottet G. H.: Equations de Navier-Stokes dans le plan avec tourbillon initial mesure.
{\it C. R. Acad. Sci. Paris S\'er. I Math.} 303, 105-108 (1986).

\bibitem{Fu}Funaki T.: A certain class of diffusion processes associated with nonlinear parabolic equations. {\it Prob. Theory and Relat. Fields}, 67(3):331-348,1984.

\bibitem{Ga-Ga}Gallagher I. and Gallay T.: Uniqueness for the two-dimensional Navier–Stokes equation with a measure as initial vorticity.
{\it Math. Ann.} 332, 287-327 (2005).


\bibitem{Gi-Mi-Os}Giga Y., Miyakawa T. and Osada H.: Two-Dimensional Naveri-Stokes flow with
measures as initial vorticity. {\it Arch. Rational Mech. Anal.}, 104, 223-250 (1988).

\bibitem{Hu-Wa}Huang X. and Wang F.Y.: Distribution dependent SDEs with singular coefficients. 
{\it Stochastic Process. Appl.} 129, no. 11, 4747-4770 (2019).

\bibitem{Ja-Wa}Jabin P.E. and Wang Z.: Quantitative estimates of propagation of chaos for stochastic systems with $W^{-1,\infty}$ 
kernels. {\it Invent. math.} 214:523-591(2018).

\bibitem{Li-Mi} Li J. and Min H.: Weak solutions of mean-field stochastic differential equations and application to zero-sum stochastic differential games. {\it SIAM Journal on Control and Optimization}, 54(3):1826-1858 (2016).

\bibitem{Ma-Be}Majda A.J. and Bertozzi A.L.: {\it Vorticity and impressible flow}. 
Cambridge Texts in Applied Mathematics, Cambridge University Press, (2002).

\bibitem{Mc} McKean H. P.: A class of Markov processes associated with nonlinear parabolic equations. {\it Proc Nat. Acad Sci USA}, 56(6):1907-1911, 1966.

\bibitem{Mi-Ve}Mishura Y.S. and Veretennikov A.Y.: Existence and uniqueness theorems for 
solutions of McKean-Vlasov stochastic equations, arXiv:1603.02212v4.

\bibitem{Ro-Zh}R\"ockner M. and Zhang X.: Well-posedness of distribution dependent SDEs with singular drifts. 
to appear in {\it Bernoulli} (2020+).

\bibitem{St-Va}Stroock D.W. and Varadhan S.S.: {\it Multidimensional diffusion processes}. Springer-Verlag, Berlin, 1979. 	

\bibitem{Sz}Sznitman A.S.: {\it Topics in propagation of chaos}. In \'Ecole d'\'Et\'e de Prob. de Saint-Flour XIX-1989, 
Vol. 1464, {\it Lect. Notes in Math.}, pages 165-251. Springer-Verlag, 1991.

\bibitem{Tr}Trevisan D.: Well-posedness of multidimensional diffusion processes with weakly differentiable coefficients. 
{\it Electron. J. Probab.\bf 21}, (2016), Paper No. 22, 41 pp.

\bibitem{Zh-Zh1}Zhang X. and Zhao G.: Singular Brownian Diffusion Processes. {\it Communications in Mathematics and Statistics,} pp.1-49, 2018. 

\bibitem{Zh-Zh2}Zhang X. and Zhao G.: Stochastic Lagrangian path for Leray solutions of 3D Navier-Stokes
  equations,  to appear in {\it Comm. in Math. Phys.} (2020+),  arXiv: 1904.04387.


\end{thebibliography}
\end{document}